\documentclass[12pt]{iopart}

\usepackage{cite}

\begin{document}

\title{On Perturbation theory for the Fokker--Planck equation}
\author{Francisco M. Fern\'{a}ndez}

\address{INIFTA (UNLP, CCT La Plata--CONICET), Divisi\'on
Qu\'imica Te\'orica, Blvd. 113 S/N,  Sucursal 4, Casilla de Correo
16, 1900 La Plata, Argentina}\ead{fernande@quimica.unlp.edu.ar}

\maketitle

\begin{abstract}
We discuss a recent application of the Modified Homotopy Perturbation Method
to the Fokker--Planck equation and show that the selected examples do not
have any connection with actual physical problems.
\end{abstract}

\section{Introduction}

In a recent paper Jafari and Aminataei\cite{JA09} proposed the application
of the homotopy perturbation method (HPM) to the Fokker--Planck equation.
More precisely, they applied a modified HPM (MHPM) to several examples and
concluded that ``we infer that our new MHPM and its coincident HPM is a
powerful tool for solving a Fokker--Planck equation.'' The purpose of this
paper is to discuss those results. In Sec.~\ref{sec:F-P} we outline the
problem. In Sec.~\ref{sec:HPM} we briefly consider the perturbation method
proposed by the authors. In Sec.~\ref{sec:Experiments} we examine each
example and its result. In Sec.~\ref{sec:convergence} we analyze the
authors' conclusions about the convergence of their method. Finally, in Sec.~%
\ref{sec:conclusions} we draw our own general conclusions.

\section{The Fokker--Planck equation}

\label{sec:F-P}

Jafari and Aminataei\cite{JA09} discussed Fokker--Planck equations of the
form
\begin{equation}
\frac{\partial u(\mathbf{x},t)}{\partial t}=-\left[ \sum_{i=1}^{N}\frac{%
\partial }{\partial x_{i}}A_{i}(\mathbf{x},t)+\sum_{i,j=1}^{N}\frac{\partial
^{2}}{\partial x_{i}\partial x_{j}}B(\mathbf{x},t)\right] u(\mathbf{x}%
,t),\;\,u(\mathbf{x},0)=f(\mathbf{x})  \label{eq:Fokk-Planck}
\end{equation}
where $\mathbf{x}=(x_{1},x_{2},\ldots ,x_{N})$. The Fokker--Planck equation
describes the time evolution of the probability density of the position of a
particle or other observables. Therefore, $u(\mathbf{x},t)$ should be an
integrable function that satisfies the boundary conditions of the problem.
In general it is not easy to solve Eq.~(\ref{eq:Fokk-Planck}) for a given
plausible initial condition $f(\mathbf{x})$.

Many textbooks discuss applications of the Fokker--Planck and diffusion
equations; one of the simplest illustrative examples is given by\cite{MQ00}
\begin{equation}
\frac{\partial C(x,t)}{\partial t}=D\frac{\partial ^{2}C(x,t)}{\partial x^{2}%
}  \label{eq:diffusion}
\end{equation}
where $C(x,t)$ is the concentration of the species and $D$ is a constant
diffusion coefficient. If we assume that $C(x,0)=C_{0}\delta (x)$, where $%
\delta (x)$ is the Kronecker delta function, then
\begin{equation}
C(x,t)=\frac{C_{0}}{2(\pi Dt)^{1/2}}\exp \left( -\frac{x^{2}}{4Dt}\right)
\label{eq:C(x,t)}
\end{equation}

\section{The homotopy story}

\label{sec:HPM}

In order to solve the Fokker--Planck equation (\ref{eq:Fokk-Planck}) Jafari
and Aminataei\cite{JA09} resorted to the HPM. Basically it consists of
introducing a parameter $p$ into the differential equation and expanding the
solution of the resulting equation $H(v,p)=0$ in a Taylor series about $p=0$%
: $v=v_{0}+v_{1}p+\ldots $. In other words, HPM is just an ordinary
perturbation approach with a fancy name. However, Jafari and Aminataei\cite
{JA09} went further and proposed a modified homotopy perturbation method
(MHPM) that they described as follows: ``We can extend $H(v,p)$ to $H(v,c,p)$%
. In the MHPM, we assume that the solution of $H(v,c,p)=0$ can be written as
a power series in $p$ and $c$; i.e. $v=c_{0}v_{0}+c_{1}v_{1}p+\ldots $,
where $c=[c_{1},c_{2},\ldots ]$. Adomian decomposition and spectral methods
are special cases of this homotopy (means MHPM).'' We clearly appreciate two
facts: first, there is no power series in $c$ and, second, if we redefine
the perturbation corrections as $\tilde{v}_{j}=c_{j}v_{j}$ then we have the
ordinary HPM because the perturbation corrections are determined by the
boundary conditions of the problem.

\section{The experiments}

\label{sec:Experiments}

Jafari and Aminataei\cite{JA09} conducted several experiments to show the
performance of this new approach. In what follows we discussed them.

\textbf{Experiment 1}:

The authors chose the trivial nonlinear differential equation $\partial
u/\partial x+u^{2}=0$, $u(0)=1$ that admits the exact solution $u(x)=1/(1+x)$%
. Notice that this equation does not have any connection whatsoever with the
Fokker--Planck one (\ref{eq:Fokk-Planck}). Therefore, it is unlikely that
this example gives us any clue about the performance of the HPM on a real
situation.

After solving the HPM equations they realized that the perturbation
corrections $v_{j}=(-x)^{j}$ are merely the terms of the Taylor expansion of
the solution about $x=0$. Of course, if one substitutes the Taylor series $%
u(x)=1+u_{1}x+u_{2}x^{2}+\ldots $ into the differential equation one easily
obtains a recurrence relation for the coefficients:
\begin{equation}
u_{n+1}=-\frac{1}{n+1}\sum_{j=0}^{n}u_{j}u_{n-j},\,n=0,1,\ldots ,\,u_{0}=1
\end{equation}
It seems that the good old Taylor approach is more efficient than the HPM
and provides a more compact expression for the solution.

\textbf{Experiment 2}:

The second example is given by
\begin{equation}
\frac{\partial u}{\partial t}=\hat{L}u,\,\hat{L}=\frac{\partial }{\partial x}%
+\frac{\partial ^{2}}{\partial x^{2}},\;u(x,0)=x  \label{eq:Exp2}
\end{equation}
where we have introduced the linear operator $\hat{L}$ in order to
facilitate the following discussion. The first question that one may ask is:
what is the physical meaning of that initial condition?. Obviously, $u(x,0)$
is unbounded and cannot be a good candidate for a probability density. It is
not difficult to guess why the authors chose it: the solution follows
straightforwardly from $u(x,t)=e^{t\hat{L}}u(x,0)=x+t$. In other words, you
cannot miss it because $\hat{L}^{n}u(x,0)=0$ for all $n>1$. It is not
surprising that the solution is also unbounded and does not resemble any
physical probability density, like, for example, Eq. (\ref{eq:C(x,t)}).

\textbf{Experiment 3}:

In this case the authors chose coefficients $A$ and $B$ that are rather
complicated and arbitrary functions of $x$ and $t$. There is no physical
justification for their choice or for the initial condition, except that the
differential equation admits the exact solution $u(x,t)=\sinh (x)e^{t}$. In
this case the authors managed to derive the Taylor series for $e^{t}$ about $%
t=0$. Once again we are in front of a solution with no physical meaning
whatsoever.

For brevity I will summarize the results of the remaining experiments.
\textbf{Experiment 4}: Taylor series about $t=0$ of $u(x,t)=xe^{t}$. \textbf{%
Experiment 5}: Taylor series about $t=0$ of $x^{2}e^{t}$. \textbf{Experiment
6}: Taylor series about $t=0$ of $(1+x)e^{t}$. We appreciate that the
authors merely chose tailor--made toy problems fabricated with the only
purpose of obtaining the known answer. Besides, their HPM always yielded the
Taylor series that one derives straightforwardly as illustrated in \textbf{%
Experiment 1}.

\section{Convergence}

\label{sec:convergence}

The authors also tried to prove the convergence of their method. However,
they simply copied the equations for the Taylor series shown in most
textbooks on calculus and did not arrive to any useful conclusion. For
example, they failed to tell the reader that the Taylor series for $1/(1+x)$
converges for $|x|<1$ or that the expansion for $e^{t}$ converges for all $t$%
.

\section{Conclusions}

\label{sec:conclusions}

I think that after the short discussion given above it is unnecessary to add
that the paper of Jafari and Aminataei\cite{JA09} should not have been
published in the first place. It may sound politically incorrect but we
think that is high time to leave euphemisms aside and discuss seriously what
is happening in today's science. There has recently been far too many such
papers published elsewhere\cite
{SNH07,CHA07,CH07a,EG07,SNH08,RDGP07,SNH08b,CH08,
M08,OA08,RAH08,SG08,ZLL08,KY09}. We will summarize the main results for the
reader's benefit. For example, Chowdhury and Hashim\cite{CH07a} applied the
HPM to obtain the Taylor series about $x=0$ of the functions $y(x)=e^{x^{2}}$%
, $y(x)=1-x^{3}/3!$, $y(x)=\sin (x)/x$, and $y(x)=x^{2}+x^{8}/72$. By means
of the same method Chowdhury et al\cite{CHA07} derived the Taylor series
about $t=0$ for the solutions of the simplest population models. Bataineh et
al\cite{SNH07} went a step further and resorted to the even more powerful
homotopy analysis method (HAM) and calculated the Taylor expansions about $%
x=0$ of the following two--variable functions: $y(x,t)=e^{x^{2}+\sin t}$, $%
y(x,t)=x^{2}+e^{x^{2}+t}$, $y(x,t)=x^{3}+e^{x^{2}-t}$, $%
y(x,t)=t^{2}+e^{x^{3}}$, $y(x,t)=-2\ln (1+tx^{2})$ and $y(x,t)=e^{-tx^{2}}$.
Zhang et al\cite{ZLL08} derived the first two terms of the Taylor expansion
about $t=0$ of $u(x,t)=-2\sec h^{2}[(x-2t)/2]$, and $u(x,t)=-(15/8)\sec
h^{2}[(x-5t/2)]$. Bataineh et al\cite{SNH08b} modified HAM to produce MHAM
and found the Taylor expansions of the functions $u(t)=t^{2}-t^{3}$ and $%
u(t)=1+t^{2}/16$. Sadighi and Ganji\cite{SG08} calculated the Taylor
expansions about $t=0$ of $u(x,t)=1+\cosh (2x)e^{-4it}$ and $%
u(x,t)=e^{3i(x+3t)}$ by means of HPM and Adomian decomposition method (ADM),
and verified that the results were exactly the same. By means of HPM Rafiq
et al\cite{RAH08} also derived polynomial functions like $y(x)=x^{4}-x^{3}$,
$y(x)=x^{2}+x^{3}$ and $y(x)=x^{2}+x^{8}/72$. \"{O}zis and Agirseven\cite
{OA08} expanded $u(x,t)=x^{2}e^{t}$, $u(x,y,t)=y^{2}\cosh t+x^{2}\sinh t$
(and other such functions) about $t=0$ by means of the HPM. Bataineh et al%
\cite{SNH08} used HAM to obtain expansions about $x=0$ for $%
w(x,t)=xe^{-t}+e^{-x}$, $w(x,t)=e^{x+t+t^{2}}$, $w(x,t)=e^{t+x^{2}}$ and $%
w(x,t)=e^{t^{2}+x^{2}}$. They thus managed to reproduce earlier ADM results.
Ko\c{c}ak and Y\i ld\i r\i m\cite{KY09} applied HPM to a 3D Green's function
for the dynamic system of anisotropic elasticity.

There are many more articles where the authors applied HPM, HAM and ADM and
produced results that anybody would easily obtain by means of a
straightforward Taylor expansion of the model differential equations. In
some cases the examples are tailor--made toy models as those discussed
above, in others the Taylor series is unsuitable for an acceptable
description of the physics of the system. In our opinion this unhappy
situation is becoming rather preoccupying, to say the least. Notice that
Frank\cite{FRANK08} criticized a previous article on the Fokker--Planck
equation in this same journal. Unfortunately, the production of such kind of
papers is inexhaustible.

The reader may find the discussion of other articles elsewhere\cite
{F07,F08b,F08c,F08d,F08e,F08f,F09a,F09b,F09c,F09d,F09e,F09f}. We recommend
the most interesting case of the predator--prey model that predicts a
negative number of rabbits\cite{F08d}.

\end{document}